\documentstyle[12pt, amsfonts]{article}
\newcommand{\ben}{\begin{enumerate}}
\newcommand{\een}{\end{enumerate}}
\newcommand{\be}{\begin{equation}}
\newcommand{\ee}{\end{equation}}
\newcommand{\bea}{\begin{eqnarray}}
\newcommand{\eea}{\end{eqnarray}}
\newcommand{\bc}{\begin{center}}
\newcommand{\ec}{\end{center}}

\newcommand{\ba}{\begin{eqnarray}}
\newcommand{\na}{\end{eqnarray}}
\newcommand{\ban}{\begin{eqnarray*}}
\newcommand{\nan}{\end{eqnarray*}}

\newcommand{\C}{{\Bbb C}} 
\newcommand{\Nr}{{\Bbb N}_r} 
\newcommand{\cP}{{\cal P}} 
\newcommand{\Tq}{{\cal T}_q({\frak g})} 
\newcommand{\cT}{{\cal T}_q({\frak g})} 
\newcommand{\Eq}{{\cal E}_q^{\frak l}} 
\newcommand{\Eg}{{\cal E}_q} 
\newcommand{\Hq}{{\cal E}_q^{\frak l}}

\newcommand{\g}{{\frak g}} 
\newcommand{\Uq}{\mbox{U}_q({\frak g})}

\newcommand{\Ul}{\mbox{U}_q({\frak l})} 
\newcommand{\Up}{\mbox{U}_q({\frak p})} 
 
\def\1{\hbox{ 1\kern-.35em\hbox{1}}}

\newcommand{\quea}{quantized universal enveloping algebra }

\newtheorem{lemma}{Lemma}
\newtheorem{proposition}{Proposition}
\newtheorem{theorem}{Theorem}
\newtheorem{definition}{Definition}

\newcommand{\cb}{\mbox{$\cal B$}}

\newcommand{\cl}{{\cal L}^2_q}

\newcommand{\ca}{{\cal A}_q({\frak g})}

\newcommand{\op}{\mbox{\scriptsize op}}

\newenvironment{proof}[1]{\begin{trivlist} \item[] {\em #1\/}: }%
{\hfill $\Box$ \end{trivlist}}

\begin{document}

\title{\normalsize{\bf COMMENTS ON THE NONCOMMUTATIVE DIFFERENTIAL 
GEOMETRY OF QUANTUM HOMOGENEOUS VECTOR BUNDLES}}
\author{R. B. Zhang\\ 
\small Department of Pure Mathematics\\
\small University of Adelaide, Adelaide, Australia\\
\small E - mail: rzhang@maths.adelaide.edu.au}
\date{}
\maketitle
\begin{abstract}
Differential calculi are obtained for 
quantum homogeneous spaces by extending Woronowicz' 
approach to the present context.  Representation theoretical 
properties of the differential calculi are investigated. 
Connections on quantum homogeneous vector bundles are 
classified  and explicitly constructed 
by using the theory of projective modules.
\end{abstract}
\vspace{1cm}
\section{\normalsize{\bf INTRODUCTION}}
The important interplay between classical differential geometry 
and the representation theory of Lie groups appears to carry 
over to the quantum setting.  
It is our broad aim to investigate this interplay.
Quantum homogeneous vector bundles provide a natural framework for 
such investigations. As shown in \cite{I}, induced 
representations of quantum groups are characterized by geometric data
on quantum homogeneous vector bundles through a generalized 
Borel-Weil theorem.
Also, the induced representations play an important role in the study 
of invariant $q$ - difference operators \cite{Dobrev}.   

This paper reports some elementary results on the differential geometry 
of quantum homogeneous vector bundles.  These bundles were introduced 
in \cite{I} following the general strategy of Connes' 
theory of noncommutative geometry \cite{Connes}.
We work in the context of Woronowicz type compact 
quantum groups \cite{groups}, thus the bundles themselves as well as their 
base spaces, the quantum homogeneous spaces,  carry natural topologies.

We will construct left covariant differential calculi on 
quantum homogeneous spaces by extending Woronowicz'
approach \cite{Woronowicz} to the present context, 
and study their representation theoretical properties.  
Connections on quantum homogeneous vector 
bundles will be classified.  We will develop an explicit construction 
of the connections by using the theory of projective modules.   

Th arrangement of the paper is as follows. Section 2 recalls 
the definition of the quantum homogeneous vector bundles.
Section 3 develops the basics of the differential geometry of 
quantum homogeneous vector bundles. The main results are 
Theorems 2 and 3, which we believe to be new. 

We should point out that there is considerable amount of work in the 
literature studying quantum principal bundles at the algebraic level
without the framework of Connes' theory \cite{Connes}.
We refer to \cite{Majid} for discussions and further references.
In comparison, quantum homogeneous vector bundles seem to be 
less studied so far.

\section{\normalsize {\bf  QUANTUM HOMOGENEOUS VECTOR BUNDLES}} 
This section recalls the definition 
of quantum homogeneous vector bundles \cite{I}. 
It also serves to fix our notation and convention. 

Let $\g$ be a finite dimensional complex simple Lie algebra 
of rank $r$, with the simple roots $\{\alpha_i \ | \ i\in\Nr\}$, 
where $\Nr=\{ 1, 2, ..., r\}$.  
We denote by $\cP$ the set of the integral elements of 
$\bigoplus_{i=1}^r{\Bbb R}\alpha_i$, 
and by  $\cP_+$ the set of the integral dominant elements. 
The Jimbo version  of the   \quea  $\Uq$ is defined to be the
unital associative algebra over $\Bbb C$,  generated by $\{ k^{\pm 1}_i, 
e_i, f_i\ | \ i\in\Nr\}$  subject to the standard 
relations \cite{Chari}.  
As is well known, $\Uq$ has the structure of a Hopf  algebra.  
We will denote the  co-multiplication by $\Delta$, 
the co-unit by $\epsilon$, and the antipode by $S$.
Throughout the paper, we assume that  $q$ is real positive.

We will consider only finite dimensional left $\Uq$ - modules 
which are direct sums of irreducible 
left modules of the following kind. If  
$W(\lambda)$ is such an irreducible module with highest weight 
$\lambda\in \cP_+$, 
then the $k_i$ eigenvalue of the highest weight vector is 
$q^{(\lambda, \ \alpha_i)/2}$ for all $i\in\Nr$.
Let $\Tq$ denote the subspace of $(\Uq)^*$, that is spanned      
by the matrix elements of the representations afforded
by the finite dimensional left $\Uq$ - modules just 
discussed. Then $\Tq$ forms a Hopf algebra.
We will also denote the co - multiplication of $\Tq$ by
$\Delta$, and the antipode by $S$.  They should not be confused
with the co - multiplication and antipode of the \quea.

Recall that $\Uq$ acquires a Hopf $\ast$-algebra structure when 
endowed with the following anti - involution, 
\ban  
e_i^\ast = f_i,&  \quad 
f_i^\ast= e_i, & \quad  
k_i^\ast=k_i.  \nan  
An important fact is that with respect to this $\ast$-operation,
every finite dimensional $\Uq$-module $W$ is unitary. 
This Hopf $\ast$-algebra structure of 
$\Uq$ induces a natural Hopf $\ast$ - algebra 
structure for $\Tq$. 
The normalized quantum Haar functional on $\Tq$ is positive definite
in the sense that $\int f^* f> 0$,  $\forall 0\ne f\in\Tq$.

The quantum group to be considered is a  
completion of $\Tq$ with an appropriate topology arising from  
the quantum Haar functional. 
Let $||\cdots||_h$ be the norm on $\Tq$ determined by
\ban
||a||_h^2 &=& \int a^* a,  \quad a\in\Tq, 
\nan
and denote by $\cl$ the Hilbert space completion of $\Tq$ 
in this norm.
Then the left regular representation of $\Tq$ can be extended to 
the completion $\cl$, yielding a 
$\ast$-representation $\phi:\Tq \to \cb(\cl)$ in 
the bounded linear operators on $\cl$.  
Let $||\cdot||$ be the operator norm on $\cb(\cl)$.  
Its pull back under $\phi$ gives rise to a $C^\ast$-semi-norm
$||\cdot ||_{\op}$ on $\Tq$.
The completion in $||\cdot ||_{\op}$  extends $\Tq$ to a
unital $C^\ast$-algebra $\ca$, which is a quantum group 
of the kind introduced by Woronowicz \cite{groups}. 

The \quea $\Uq$ acts on $\cT$ naturally. 
We will make extensive use of two actions, `$\circ$'
and `$\cdot$', which are respectively defined by 
\ban
x\circ f&=&\sum_{(f)} f_{(1)} \ \langle f_{(2)}, \ x\rangle,\\ 
x\cdot f &=&\sum_{(f)} \langle f_{(1)}, \ S^{-1}(x) \rangle f_{(2)}, 
\quad x\in\Uq, \ f\in\cT.   \label{circ}
\nan
Note that $\circ$
corresponds to the right translation in the classical theory
of Lie groups, while $\cdot$ corresponds to the
left translation.  The actions commute.

Using these two actions, we introduce the 
quantum analog of the algebra of smooth functions over a Lie group:
\begin{definition}
\ba
{\cal E}_q:= \{ a\in\ca |  x\cdot a,\ x\circ a\in \ca,\
                \ |a(x)|<\infty, \ \forall x\in\Uq\}.
\na
\end{definition}
Clearly, ${\cal E}_q$ forms a subalgebra of $\ca$.

For any subset $\bf\Theta$ of $\Nr$, we
introduce the following two sets of elements of $\Uq$:
\ban
{\cal S}_l=\{ k_i^{\pm 1}, i\in \Nr; \
\ e_j, \ f_j, \ j \in {\bf\Theta}\};
\quad \quad
{\cal S}_{p}= {\cal S}_l \cup \{ e_j,
      j \in \Nr\backslash {\bf\Theta}\}.
\nan
${\cal S}_l$ generates a reductive  quantum subalgebra $\Ul$,
while ${\cal S}_{p}$  generates a parabolic quantum  
subalgebra $\Up$ of $\Uq$. 

Corresponding to each $\Ul$, we introduce a quantum homogeneous 
space, on which the algebra of `smooth' functions is given by  
\begin{definition}    
\ba 
\Eq&:=&\left\{ f\in\Eg \ | \ x\circ f =\epsilon(x) f, 
\quad \forall x\in\Ul \right\}. \label{space} 
\na
\end{definition}
Note that $\Eq$ is an infinite dimensional subalgebra of $\Eg$.
Furthermore, $\Eq$ forms a left $\Eg$ co-module under the 
co-multiplication of $\Eg$. 

Let $V$ be a finite dimensional left $\Ul$-module.
We extend the actions $\circ$ and $\cdot$ of $\Uq$ on $\Eg$ trivially 
to actions on $V\otimes_\C \Eg$: for any 
$\zeta=\sum_r v_r\otimes f_r$  $\in V\otimes \Eg$   
\ban 
x\circ\zeta = \sum_r  v_r \otimes x\circ f_r,&   
x\cdot\zeta = \sum_r  v_r \otimes x\cdot f_r, & \quad   x\in\Uq. 
\nan  
The quantum homogeneous vector bundle induced from the $\Ul$ - module $V$ 
is specified by its space of sections, which is defined by    
\begin{definition}
\ba 
\Hq(V)&:=& \left\{ \zeta\in V\otimes_\C \Eg\ 
| \ x\circ\zeta = (S(x)\otimes id_{\Eg}) \zeta, \ \forall x\in\Ul\right\}.
\na  
\end{definition}
$\Hq(V)$ is infinite dimensional if any weight of $V$ belongs to $\cP$, 
and zero otherwise. $\Ul$ being reductive, $V$ can be reduced 
into a direct sum of irreducible $\Ul$ - modules.  An irreducible 
module does not yield any nonzero section in $\Hq(V)$  if any 
of its weights does not belong to $\cP$.   Therefore, we can assume 
that all the weights of $V$ are integral. 

We define a left and a right action of $\Eq$ on $\Hq(V)$. 
For all $a\in\Eq$, and $\zeta=\sum_r v_r \otimes f_r\in\Hq(V)$, 
\ban
a\otimes \zeta\mapsto a\zeta&=&\sum_r  v_r \otimes a f_r, \\ 
\zeta \otimes a \mapsto \zeta a&=&\sum_r v_r \otimes f_r a. 
\nan
It was shown in \cite{I} that 
\begin{theorem}\label{proj}
$\Hq(V)$ furnishes  a two - sided $\Eq$ - module, which 
is projective and of finite type  
both as a left and a right module 
(but not as a two - sided module.). 
\end{theorem}

Therefore, our definition of quantum homogeneous vector 
bundles is consistent with the general definition of 
noncommutative vector bundles in Connes' theory \cite{Connes}. 
As we will see in the next section, 
the theorem is also of major importance for developing  
noncommutative differential geometry on quantum homogeneous 
vector bundles.
In fact, it follows from a theorem of \cite{Cuntz} that a connection
exists on the quantum homogeneous vector bundle if and only if 
$\Hq(V)$ is projective as a right (or
left) $\Eq$ - module.

A further remark is that when $V$ is a finite dimensional 
irreducible $\Up$ - module, which naturally restricts to an
irreducible $\Ul$ - module, 
the quantum analog of the space of `holomorphic' sections 
\ban
{\cal O}_q(V)&=&\left\{\zeta\in\Hq(V)\, |\, 
p\circ\zeta = (S(p)\otimes \mbox{id}_{\Eg})\zeta, \ \forall p\in\Up\right\},
\nan 
forms an irreducible or zero $\Uq$ - module under the action $\cdot$ \cite{I}.  

\section{\normalsize {\bf NONCOMMUTATIVE DIFFERENTIAL GEOMETRY}}
\subsection{Differential calculi on quantum homogeneous spaces} 
Let 
\bea 
d: \Eg&\rightarrow& \Gamma  \label{first}
\eea 
be a bi-covariant first order 
differential calculus of $\Eg$ \cite{Woronowicz}. 
$\Gamma$ is by definition a bi-covariant two-sided $\Eg$ - module. 
Recall that the compatibility of
the left and right $\Eg$ - co-module structures of $\Gamma$ 
guarantees that the left co-invariant subspace $\Gamma_{inv}$ 
of $\Gamma$ has a 
right $\Eg$ - co-module structure, which in turn induces a   
left $\Uq$ - module structure. 
We assume that $\Gamma_{inv}$ is of dimension $K<\infty$  as a 
complex vector space.
Now $\Gamma$ is a free left $\Eg$ - module, and 
any basis of $\Gamma_{inv}$ gives rise to a basis
of $\Gamma$ over $\Eg$. 
Denote by $\pi$ the $\Uq$ - representation associated with  
$\Gamma_{inv}$ in the basis.  Let $\sigma =$  
$P(\pi\otimes \pi)R$ with $P$ the flipping map and $R$ the universal 
$R$ - matrix of $\Uq$.     
Then $\sigma$ can be decomposed in a unique fashion into the form 
$\sigma_+ - \sigma_-$, where $\sigma_\pm$ are $\Uq$ - module 
maps, which satisfy the relations $\sigma_+ \sigma_- $
$=\sigma_- \sigma_+ =0$, and respectively reduce to the 
projections onto the symmetric 
and skew symmetric tensors in the classical limit.    
We extend $\sigma_-$ to an endomorphism of $\Eg\otimes_\C 
\Gamma_{inv}\otimes_\C \Gamma_{inv}$ by requiring it to act on 
the space $\Eg$ trivially. This defines a map on 
$\Gamma\otimes_{\Eg}\Gamma$. 

Denote by $\Gamma^\otimes$ the tensor algebra 
\ban
    \Gamma^\otimes &=& \bigoplus_{n\ge 0} \Gamma^{\otimes n},\\
         \Gamma^{\otimes n}&=&\underbrace{\Gamma\otimes_{\Eg}
   \Gamma\otimes_{\Eg} ... \otimes_{\Eg} \Gamma}_n.    
\nan 
Let $J$ be the two-sided ideal of $\Gamma^\otimes$ generated by 
$ker\sigma_-$. We define 
\ban
\Omega(\Eg)&=&\Gamma^\otimes/J, 
\nan
and denote by $\tau$ the canonical epimorphism $\Gamma^\otimes
\rightarrow \Omega(\Eg)$. Note that $\Omega(\Eg)$ is graded. 
Some elementary considerations of $\Uq$ -  representations  
lead to the conclusion that $\tau(\Gamma^{\otimes n})=0$, $\forall n>K$.
Therefore, $\Omega(\Eg)=\oplus_{n=0}^K \Omega^n(\Eg)$. 
It is known that there exists a unique bi-covariant differential 
calculus \cite{Woronowicz}  
\bea 
d: \Omega(\Eg) &\rightarrow & \Omega(\Eg), \label{differential}
\eea  
such that $d: \Eg\rightarrow \Omega^1(\Eg)$ coincides 
with the given first order differential calculus.
The study of bi-covariant differential calculi on quantum groups 
has been actively pursued in recent years, and a good 
understanding of them has been gained. We refer to \cite{Schmudgen} 
for details and further references, and in particular, 
the third  paper of \cite{Schmudgen}  for a discussion of 
the problem of extending a first order differential calculus 
on a quantum group to a higher order one.

Recall that we defined the algebra of functions $\Eq$ 
on a quantum homogeneous space as a subalgebra of $\Eg$. 
Thus the first order differential calculus (\ref{first}) on 
the quantum group can be naturally restricted to  the quantum 
homogeneous space $\Eq$.  We will still use the 
same notation $d$ to denote the restriction of $d$ to $\Eq$.
Set 
\bea 
\Omega^1(\Eq)&=&\{\sum a_i d b_i\, | \, a_i, b_i\in\Eq\}.
\eea  
Then 
\begin{proposition}
The following linear map gives rise to a 
left $\Eg$ - covariant first order differential calculus on $\Eq$:
\bea d: \Eq\rightarrow \Omega^1(\Eq).\label{d} \eea 
\end{proposition}
Note that $\Omega^1(\Eq)$ forms a two-sided 
$\Eq$ - module,  and also a left $\Eg$ -co-module. 
It follows from the compatibility of left $\Eg$ - co-module
and left $\Eg$ - module structures of $\Eq$ 
that $\Omega^1(\Eq)$ is left covariant with respect to $\Eg$.  

Introduce the tensor algebra 
\bea
\Omega^\otimes(\Eq)&=&\bigoplus_{n\ge 0} (\Omega^1(\Eq))^{\otimes n},
\nonumber\\
(\Omega^1(\Eq))^{\otimes n}&=& \underbrace{\Omega^1(\Eq)\otimes_{\Eq} 
\Omega^1(\Eq)\otimes_{\Eq} ... \otimes_{\Eq} 
\Omega^1(\Eq)}_n, 
\eea 
and define 
\bea 
\Omega(\Eq)&=&\bigoplus \Omega^n(\Eq),\nonumber\\
\Omega^n(\Eq)&=&\tau((\Omega^1(\Eq))^{\otimes n}).
\eea 
Then we have 
\begin{theorem}
The following linear map defines a 
left $\Eg$ - covariant differential calculus on $\Eq$:
\bea d: \Omega(\Eq) \rightarrow \Omega(\Eq), \label{calculus}\eea
where $d$ denotes the restriction of (\ref{differential}) to 
$\Omega(\Eq)$. 
\end{theorem}
A moment's reflection will convince ourselves that 
it is indeed true that $d(\Omega(\Eq))$ is a subset of  
$\Omega(\Eq)$.   Needless to say, this differential calculus 
has the standard properties, namely, 
\begin{description}
\item i). $d: \Eq\rightarrow \Omega^1(\Eq)$ coincides with (\ref{d});
\item ii). $d: \Omega^n(\Eq) \rightarrow \Omega^{n+1}(\Eq)$;
\item iii). $d^2=0$;
\item iv). $d(\theta \omega) = (d\theta) \omega +
           (-1)^n \theta (d\omega)$, $\forall$  $\theta\in \Omega^n(\Eq)$, 
$\omega\in \Omega(\Eq)$. 
\end{description}

The left $\Eg$ - module structure of $\Omega(\Eq)$ 
induces a left $\Uq$ - module structure. 
Let $\delta_L: \Omega(\Eq)$ $\rightarrow \Eg\otimes_\C \Omega(\Eq)$ 
denote the left $\Eg$ - co-action.  
Then the left $\Uq$ action is given by  
\ban 
\cdot :  \Uq\otimes \Omega(\Eq)&\rightarrow& \Omega(\Eq),\\
x\otimes\omega &\mapsto& x\cdot\omega= \delta_L(\omega)(S^{-1}(x)).
\nan 
The right covariance of the differential calculus on 
$\Eg$ is no longer preserved when restricted to the quantum 
homogeneous space. However, a residue of it implies that 
$\Omega(\Eq)$ forms a left $\Ul$ - module:
\ban
\circ : \Ul \otimes \Omega(\Eq) &\rightarrow& \Omega(\Eq),\\  
         p\otimes \omega &\mapsto& p\circ \omega,  
\nan  
where $\circ$ is defined on $\Omega^1(\Eq)$ by 
$ p\circ(\sum_i a_i d b_i)= \sum_i\sum_{(p)} (p_{(1)}\circ a_i ) 
d(p_{(2)}\circ b_i )$ and extended to the entire $\Omega(\Eq)$
in the obvious way.  Clearly, $p\circ\omega =\epsilon(p)\omega$, 
$\forall p\in\Ul$, $\omega\in\Omega(\Eq)$.  

The bi-covariance of the differential calculus (\ref{differential}) 
on the quantum group implies that the differential operator 
(\ref{d}) is a module homomorphism with respect to the afore 
discussed $\Uq$ and $\Ul$ module structures of $\Omega(\Eq)$. 
That is
\begin{lemma}
Given any $\omega\in\Omega(\Eq)$, 
\ban 
x\cdot(d\omega) = d (x\cdot\omega),&
p\circ(d\omega)= d(p\circ\omega), & \forall x\in\Uq, \ p\in\Ul.
\nan 
\end{lemma}

\subsection{Connections on quantum homogeneous vector bundles} 
In order to construct connections on a quantum homogeneous vector 
bundle, we need to better understand the module structures of
its space of sections with respect to the algebra of functions 
on the quantum homogeneous space. Recall that
given a finite dimensional left $\Ul$ - module $V$ with integral weights,
there always exists a finite dimensional left $\Ul$ - module $V^\bot$
such that $V\oplus V^\bot$ is isomorphic to
the restriction of some left $\Uq$ - module $W$. 
We emphasize that there exists
a straightforward procedure for constructing a required $W$
(though $W$ may not be unique).  
Let $\Im_o: V\rightarrow W$
be the embedding, and $\wp_o: W\rightarrow V$ be the canonical
projection.  Clearly $\wp_o \Im_o = \mbox{id}_V$.  Define
\ban
\wp: W\otimes_\C \Eq&\rightarrow&\Hq(V),
\nan
by the composition of maps
\bea
W\otimes_\C \Eq \stackrel{\delta\otimes id}{\longrightarrow}
W\otimes_\C\Eg\otimes_\C\Eq
\stackrel{ id\otimes M(S\otimes id)}{\longrightarrow}
W\otimes_\C\Eg\stackrel{\wp_o\otimes id}{\longrightarrow}\Hq(V),
\eea
where $M$ is the multiplication of $\Eg$, and
$\delta: W \rightarrow W\otimes_{\C}\Eg$ is the
$\Eg$ co-action on $W$ dualizing the left $\Uq$ action
of $W$. Explicitly, $\delta$ is defined by     
\be
\delta(w)(x) = x w, \quad  \forall w\in W,\ x\in\Uq.
\label{comodule}
\ee
Similarly, we define
\ban
\Im: \Hq(V)&\rightarrow& W\otimes_\C \Eq,
\nan
by
\bea
\Hq(V)\stackrel{\Im_o\otimes id}{\longrightarrow}
W\otimes_\C\Eg\stackrel{\delta\otimes id}{\longrightarrow}
W\otimes_\C\Eg\otimes_\C\Eg \stackrel{ id\otimes M}{\longrightarrow}
W\otimes_\C\Eq.
\eea
Then one can easily show that
\begin{proposition}
$\wp$ and $\Im$ are well defined right $\Eq$ - module 
homomorphisms.  $\wp$ is surjective,  $\Im$ is injective, 
and $\wp \Im=\mbox{id}_{\Hq(V)}$.
\end{proposition}

With the help of these maps, we can now define 
a linear map 
\bea
\partial: \Hq(V)&\rightarrow& \Hq(V)\otimes_{\Eq} \Omega^1(\Eq),
\label{partial}
\eea
which sends sections of a quantum homogeneous vector bundle 
to the tensor product of the space of sections of the bundle 
with the  space of one - forms on the base space.
It is defined by the composition of the following maps:
\bea
\Hq(V)\stackrel{\Im}{\longrightarrow}
W\otimes_{\C}\Eq \stackrel{id \otimes d}{\longrightarrow}
W\otimes_{\C}\Omega^1(\Eq)\stackrel{\wp \otimes id}{\longrightarrow}
\Hq(V)\otimes_{\Eq} \Omega(\Eq). 
\eea
An important property of $\partial$ is that for all
$\zeta\in\Hq(V)$ and $a\in\Eq$,  $\partial(\zeta a) 
=\partial(\zeta) a + \zeta d a$. 
This allows us to extend $\partial$ to a linear map 
\bea
\nabla_o:
\Hq(V)\otimes_{\Eq} \Omega(\Eq) &\rightarrow&
\Hq(V)\otimes_{\Eq} \Omega(\Eq), 
\eea
by requiring that 
\ban 
\nabla_o(\zeta \omega)&=& \partial(\zeta) \omega + \zeta d\omega,
\quad \forall \zeta\in\Hq(V), \ \omega\in\Omega(\Eq).
\nan
Then
\begin{theorem}
$\nabla_o$ defines a connection on the quantum homogeneous vector 
bundle $\Hq(V)$.  Furthermore, every connection on $\Hq(V)$ is of 
the form $\nabla_o+A$ for some $A\in End_{\Eq}(\Hq(V))$.
\end{theorem}
\begin{proof}{\em Proof}: 
The Theorem, being rather obvious, hardly requires any proof. We 
nevertheless make some explanatory remarks below. 
Recall that a connection on the bundle $\Hq(V)$ is by definition 
a linear map 
\ban 
{\nabla}: 
\Hq(V)\otimes_{\Eq} \Omega(\Eq) &\rightarrow&
\Hq(V)\otimes_{\Eq} \Omega(\Eq),  
\nan
such that $\nabla(\Hq(V)\otimes_{\Eq} \Omega^n(\Eq))$
$\subset \Hq(V)\otimes_{\Eq} \Omega^{n+1}(\Eq)$, and 
\ban
\nabla(\psi\omega)&=&\nabla(\psi) \omega 
+(-1)^n \psi d\omega, 
\quad \forall \psi\in\Hq(V)\otimes_{\Eq} \Omega^n(\Eq),
\ \omega\in\Omega(\Eq).
\nan 
The map $\nabla_o$ is designed to have these properties, 
thus yielding a connection for $\Hq(V)$.

If $\nabla$ is also a connection, then 
\ban 
(\nabla-\nabla_o)(\psi\omega)&=& (\nabla-\nabla_o)(\psi)\omega.
\nan
This in particular implies that 
$\nabla-\nabla_o\in End_{\Eq}(\Hq(V))$,  
where $End_{\Eq}(\Hq(V))$ denotes the space of the right $\Eq$ -
linear maps $\Hq(V)\rightarrow \Hq(V)$. 
\end{proof}

The construction of the connections can be made entirely explicit.  
To do that we only need to consider the distinguished connection 
$\nabla_o$. 
Note that we can always choose appropriate bases 
for $W$ and $V$, respectively denoted by,   
\ban
\{w_\alpha\, |\, \alpha=1, 2, ..., dimW\}, &\quad&  
\{ v_i\, |\, i=1, 2, ..., dimV\},
\nan
such that 
\ban
\wp_o(w_\alpha)&=&\delta_{i\alpha} v_i,\\
\Im_o(v_i)&=&w_i, \quad 1\le i\le dimV,\ \ 1\le \alpha\le dimW.
\nan 
Denote by $t$ the $\Uq$ - representation afforded by 
$W$ in the given basis, and by $t_{\alpha \beta}$ the matrix 
elements.  Then for any $p\in\Ul$, $p w_i = \sum_{j=1}^{dimV} 
t_{j i}(p) w_j$. Also, $\delta(w_\alpha)=\sum_{\beta=1}^{dimW} 
w_\beta\otimes t_{\beta \alpha}$. 
The following $dimW$ sections, 
\ban 
\zeta_\alpha &= \wp(w_\alpha\otimes \epsilon)
             &= \sum_{i=1}^{dimV} v_i \otimes S(t_{i \alpha}), 
\nan 
which may not be linearly independent in general,   
generate the entire right $\Eq$ - module $\Hq(V)$, i.e., 
every section of $\Hq(V)$ can be expressed in the form 
\ban 
\zeta &=&\sum_\alpha \zeta_\alpha a_\alpha, 
\quad a_\alpha\in\Eq.
\nan    
It is not difficult to work out that 
\ban
\partial(\zeta)&=& \sum_{\beta,\, \alpha,\,  j}
\zeta_\beta d(t_{\beta j} S(t_{j \alpha}) a_{\alpha}).
\nan 
Now every $\phi\in\Hq(V)\otimes_{\Eq}\Omega(\Eq)$ is of the 
form  
\ban 
\phi&=&\sum_{\alpha} \zeta_{\alpha}\omega_{\alpha},
\quad \omega_{\alpha}\in\Omega(\Eq).
\nan
We have 
\ban 
\nabla_o(\phi)&=&\sum_{\beta,\,\alpha,\,  j} 
\zeta_\beta d(t_{\beta j} S(t_{j \alpha})) \omega_{\alpha}
+\sum_{\alpha}\zeta_\alpha d\omega_{\alpha}.
\nan 
  
Before closing, we mention that  
given any connection $\nabla$, 
$\nabla^2$ is right $\Omega(\Eq)$ - linear. The curvature $F$  
associated with this connection is defined to be the restriction 
of $\nabla^2$ to $\Hq(V)$, which satisfies the Bianchi identity. 
This of course is standard.

\end{document}